\RequirePackage{fix-cm}
\documentclass[smallextended]{svjour3}       
\smartqed  
\usepackage{amssymb,epsfig,multirow}
\usepackage{graphicx}
\usepackage{float}
 \usepackage{color}
\usepackage{enumerate}
\usepackage{hyperref}

\newtheorem{thm}{Theorem}

\begin{document}
	
\title{A  new Barzilai-Borwein steplength from the viewpoint of total least squares
\thanks{This research was supported by the Beijing Natural Science Foundation under grant Z180005, and the National Natural Science Foundation of China under grants 11822103, 11571029, and 11771056.}
}

\titlerunning{A  new  Barzilai-Borwein steplength}        

\author{Shiru Li \and Yong  Xia
}


\institute{ S. Li \and    Y.  Xia    \at
	 LMIB of the Ministry of Education, School of Mathematical Sciences, Beihang University, Beijing, 100191, P. R. China
	\email{lishiru@buaa.edu.cn  (S. Li); yxia@buaa.edu.cn (Y. Xia, corresponding author)}
}

\date{Received: date / Accepted: date}

\maketitle

\begin{abstract}
Barzilai-Borwein (BB) steplength is a popular choice in gradient descent method. By observing that the two existing BB steplengths correspond to the ordinary and the data least squares, respectively, we employ the third kind of least squares, the total least squares, to create a new BB steplength, which is shown to lie between the two existing BB steplengths.


\end{abstract}

\section{Introduction}
The classical iterative formulation of the gradient descent method for unconstrained minimization reads as
\begin{eqnarray*}
	x_{k+1}=x_{k}-\alpha_k g_k, ~k=0,1,\ldots,
\end{eqnarray*}
where $g_k$ is the gradient of the objective function at the $k$-th iteration point $x_k$, and $\alpha_k$ is a steplength with many choices.
BB steplength, proposed by Barzilai and Borwein \cite{barzilai1988two}, is one of the popular choices of $\alpha_k$. The idea is to approximate the quasi-Newton matrix
by $\alpha_k I$ (where $I$ is an identity matrix). It  leads to the simplified quasi-Newton equations and their inverse version as $s_k\approx\alpha y_k$ and $\frac{1}{\alpha}s_k\approx y_k$, respectively, where $y_k=g_{k+1}-g_k$ and $s_k=x_{k+1}-x_k$.  
Due to the over-determination, the ordinary least squares is introduced to solve the simplified equations: \begin{eqnarray}
	&&\alpha^{BB1}_{k+1}=\arg\min\limits_{\alpha} \|s_k-\alpha y_k\|^2=\frac{s^T_k y_k}{y^T_k y_k},\label{b1}\\
	&&\alpha^{BB2}_{k+1}=\frac{1}{\arg\min\limits_{\beta} \|\beta s_k-y_k\|^2}=\frac{s^T_k s_k}{s^T_k y_k},\label{b2}
\end{eqnarray}
where $\|\cdot\|$ is the standard Euclidean norm.

BB gradient method achieved many successes in practice, for
recent applications, we refer to \cite{abubakar2020barzilai,liang2019barzilai,yu2020minibatch}.
Though in general BB method is not convergent \cite{fletcher2005barzilai}, many studies focused on convergence analysis on strongly convex quadratic functions \cite{barzilai1988two,dai2013new,dai2002r,dai2005asymptotic,raydan1993barzilai,Yuan1993},
and variant BB method
 \cite{dai2006cyclic,raydan1997barzilai}.

In this note, we take a new look at BB steplengths from the viewpoint of different kinds of least squares.
As we can observe that $\alpha^{BB1}$ (\ref{b1}) and $\alpha^{BB2}$ (\ref{b2}) correspond to the ordinary and the data least squares, respectively, as presented in Section 2, it implies from the total least squares a new BB steplength formula, which is located between  $\alpha^{BB1}$ and $\alpha^{BB2}$. We use an example to show the good balance of the newly proposed BB steplength. Conclusions are made in Section 3.

\section{Total least squares and the third BB steplength}
The main contribution of this section is to derive a new BB steplength from the perspective of total least squares.
\subsection{Least squares}
We briefly review three kinds of least squares in this subsection.
The goal is to solve an over-determined linear system $Ax\approx b$, where $A\in \mathbb{R}^{m\times n}$ is the data matrix, and $b\in \mathbb{R}^m$ is the observation vector.

If we assume that only $b$ contains a noise $r$, solving
\begin{eqnarray*}
\min\limits_{r,x} \{\|r\|^2:Ax=b+r\}=\min\limits_x \|Ax-b\|^2,
\end{eqnarray*}
yields the well-known ordinary least squares problem.

The less popular data least squares problem \cite{degroat1993data}  corresponds the case where only $A$ is noised by $E$:
\begin{eqnarray*}
\min\limits_{E,x} \{\|E\|^2_F:(A+E)x=b\}=\min\limits_x \frac{\|Ax-b\|^2}{\|x\|^2},
\end{eqnarray*}
where $\|\cdot\|_F$ is the Frobenius norm.

In case that both $A$ and $b$ are noised, the total least squares problem reads as follows:
\begin{eqnarray*}
\min\limits_{E,r,x}\{\|E\|^2_F+\|r\|^2:
(A+E)x=b+r\}=\min\limits_{x}\frac{\|Ax-b\|^2}{\|x\|^2+1},
\end{eqnarray*}
which was firstly proposed in \cite{golub1980analysis}.

\subsection{A new formula and its property}

First, by replacing the ordinary least squares in (\ref{b1})-(\ref{b2}) with the data least squares, we obtain $\alpha^{BB2}$ and $\alpha^{BB1}$, respectively.

Now we apply the total least squares to solve
the simplified quasi-Newton equations $s_k\approx\alpha y_k$:
\[
	\min\limits_\alpha q(\alpha)=\frac{||\alpha y_k-s_k||^2}{\alpha^2+1},
\]
which gives the solution
\begin{equation}
	\alpha^{BB3}_{k+1}=\frac{s^T_k s_k-y^T_k y_k+\sqrt{(y^T_k y_k-s^T_k s_k)^2+4(s^T_k y_k)^2}}{2s^T_k y_k}. \label{b3}
\end{equation}

Moreover, if we apply the total least squares to solve the inverse quasi-Newton equations $\beta s_k\approx y_k$ with respect to $\beta:=1/\alpha$, we can obtain the same formula (\ref{b3}).

We can reformulate $\alpha^{BB3}$ (\ref{b3}) as a function in terms of $\alpha^{BB1}$ and $\alpha^{BB2}$. It reveals how
$\alpha^{BB3}$ keeps a balance between $\alpha^{BB1}$ and $\alpha^{BB2}$. We list  these observations in the following and omit the trivial proof.
\begin{thm} \label{thm}
Suppose $\alpha^{BB1}_{k+1}>0$, then we have
\[ \alpha^{BB3}_{k+1}=\frac{\alpha^{BB2}_{k+1}-\frac{1}{\alpha^{BB1}_{k+1}}
+\sqrt{\left(\frac{1}{\alpha^{BB1}_{k+1}}-\alpha^{BB2}_{k+1}\right)^2+4}}{2}.
\]	
Moreover, it holds that
\begin{eqnarray} &&\alpha^{BB1}_{k+1}\le\alpha^{BB3}_{k+1}\le\alpha^{BB2}_{k+1},
\nonumber\\
&&\lim\limits_{\alpha^{BB1}_{k+1}\to\infty} \frac{\alpha^{BB3}_{k+1}}{\alpha^{BB2}_{k+1}}=1,~\lim\limits_{\alpha^{BB2}_{k+1}\to 0} \frac{\alpha^{BB3}_{k+1}}{\alpha^{BB1}_{k+1}}=1.	\nonumber
\end{eqnarray}
\end{thm}

Roughly speaking, BB method with $\alpha^{BB1}$ converges more robustly and hence less fast  than that with $\alpha^{BB2}$.
Theorem \ref{thm} suggests that $\alpha^{BB3}$ seems to be a balance between $\alpha^{BB1}$ and $\alpha^{BB2}$. If both $\alpha^{BB1}$ and $\alpha^{BB2}$ are small (or large) enough, $\alpha^{BB3}$ automatically approaches to the smaller (or larger) one.

We numerically show the benefit of the balance by a classical example.
Consider the two-dimensional Rosenbrock function \cite{more1981testing}
\begin{eqnarray*}
	f(x)=100(x_2-x^2_1)^2+(1-x_1)^2,~(x_1^0,x_2^0)=(-1.2,1).
\end{eqnarray*}
Starting from the same initial points $(x_1^1,x_2^1)=(x_1^0,x_2^0)$,
we independently run the three BB methods with $\alpha^{BB1}$, $\alpha^{BB2}$, and $\alpha^{BB3}$, respectively. We set the stop criterion as $\|(x^{k}_1,x^{k}_2)-(x^*_1,x^*_2)\|\le\epsilon$ together with a maximum iteration number $5000$, where $(x^*_1,x^*_2)=(1,1)$ is the minimizer of $f(x)$.
We report in Table \ref{tab:1} the iteration numbers with different setting of $\epsilon$, where ``--'' stands for the situation that the maximum iteration number is reached.

\begin{table}[h]\centering
	\caption{Comparison of iteration numbers in minimizing the planar Rosenbrock function.} \label{tab:1}
	\begin{center}
		\begin{tabular}{cccc} \hline
			$\epsilon$  &BB1  &BB2  &BB3      \\
			\hline
			$10^{-1}$ & 154& --& 32\\
			$10^{-2}$& 160& --& 38\\
			$10^{-4}$& 166& --& 44\\
			$10^{-8}$ & 172& --& 46\\
			\hline			
		\end{tabular}
	\end{center}
\end{table}

\section{Conclusions}
From the perspective of least squares, we show the two existing BB steplengths correspond to the ordinary and data least squares, respectively. Then, based on the third one, the total least squares,
we propose a new BB steplength in this note. We prove that it lies (and hence keeps a balance) between the two existing BB steplengths. Future studies include more convergence analysis and variants of BB methods based on the new steplength.

\bibliographystyle{plain}
\bibliography{2OL2}

\end{document}